\documentclass[12pt]{amsart}

\usepackage{latexsym}
\usepackage{amssymb}
\usepackage[cp850]{inputenc}
\usepackage{epsfig}
\usepackage{psfrag}
\usepackage{amsthm}
\usepackage{amscd}
\usepackage{amsmath}
\usepackage{amsfonts}
\usepackage{graphics}
\usepackage[all]{xy}

\usepackage{fancyhdr}
\newtheorem{sat}{Theorem}[section]			\newtheorem{lem}[sat]{Lemma}
\newtheorem{kor}[sat]{Corollary}			\newtheorem{prop}[sat]{Proposition}
				
\newtheorem*{defi*}{Definition}	\newtheorem*{bei*}{Example}
\newtheorem*{sat*}{Theorem}			\newtheorem*{kor*}{Corollary}
\newtheorem*{rmk*}{Remark}

\let\ssection=\section
\renewcommand{\section}{\setcounter{equation}{0}\ssection}

\newtheorem*{namedtheorem}{\theoremname}
\newcommand{\theoremname}{testing}
\newenvironment{named}[1]{\renewcommand{\theoremname}{#1}\begin{namedtheorem}}{\end{namedtheorem}}

\theoremstyle{remark}

			\newcommand{\BH}{\mathbb H}
\newcommand{\BR}{\mathbb R}			
			\newcommand{\BQ}{\mathbb Q}
\newcommand{\BS}{\mathbb S}			\newcommand{\BZ}{\mathbb Z}

		\newcommand{\CB}{\mathcal B}

\newcommand{\CS}{\mathcal S}		
		
		\newcommand{\CX}{\mathcal X}
\newcommand{\CY}{\mathcal Y}

\newcommand{\D}{\partial}

\DeclareMathOperator{\SL}{SL}		

\DeclareMathOperator{\GL}{GL}		
\DeclareMathOperator{\Id}{Id}		

\DeclareMathOperator{\Nil}{Nil}

\DeclareMathOperator{\SO}{SO}

\newcommand{\bs}{\backslash}
\DeclareMathOperator{\syst}{syst}
\DeclareMathOperator{\vcdim}{vcdim}
\DeclareMathOperator{\Span}{Span}

\begin{document}

\title[]{Minimality of the well-rounded retract}

\author{Alexandra Pettet \& Juan Souto}
\thanks{\tiny Research partially supported by the NSF grant 0706878.}

\begin{abstract}
We prove that the well-rounded retract of $\SO_n\bs\SL_n\BR$ is a minimal $\SL_n\BZ$-invariant spine.
\end{abstract}

\maketitle

\section{Introduction}

In this note we are interested in a certain $\SL_n\BZ$-invariant deformation retract of the symmetric space $S_n=\SO_n\bs\SL_n\BR$. To every element $A\in\SL_n\BR$ one can associate the lattice $A\BZ^n$ in $\BR^n$. The element $A$ is {\em well-rounded} if the set of shortest non-zero vectors of the lattice $A\BZ^n$ generate $\BR^n$ as a real vector space. This property is invariant under the left action of $\SO_n$ and hence there is no ambiguity in saying that an element in $S_n$ is well-rounded. The subset $\CX$ of $S_n$ consisting of well-rounded elements is homeomorphic to an $\frac{n(n-1)}2$-dimensional CW-complex and the right action of $\SL_n\BZ$ on $S_n$ induces a cocompact action on $\CX$. Observe that if $n=2$ then $\CX$ is the dual to the Farey tesselation of $S_2=\BH^2$ and hence homeomorphic to the Bass-Serre tree of $\SL_2\BZ$. For larger $n$, the set $\CX$ does not have such a simple description, but Soul\'e \cite{Soule} in the case of $n=3$ and Ash \cite{Ash} in general proved that $\CX$ is a deformation retract of $S_n$ and hence contractible. This is why the subset $\CX$ is known as the {\em well-rounded retract} of $S_n$. Our goal is to show that $\CX$ is a minimal $\SL_n\BZ$-invariant spine of $S_n$.

\begin{defi*}
Let $\Gamma$ be a group acting discretely on a contractible space $S$. We say that a closed subset $X$ of $S$ is a {\em minimal $\Gamma$-invariant spine} if it is $\Gamma$-invariant and contractible and does not properly contain any closed set with these properties.
\end{defi*}

\noindent We prove:

\begin{sat}\label{main}
The well-rounded retract $\CX$ is a minimal $\SL_n\BZ$-invariant spine of the symmetric space $S_n=\SO_n\bs\SL_n\BR$.
\end{sat}

It has long been known that the well-rounded retract does not contain any smaller dimensional $\SL_n\BZ$-invariant spines. This follows namely from the fact due to Borel-Serre \cite{Borel-Serre} that the group $\SL_n\BZ$ has virtual cohomological dimension 
$$\vcdim(\SL_n\BZ)=\frac{n(n-1)}2=\dim\CX$$
In order to appreciate the difference between this statement and the claim of Theorem \ref{main} it should be observed that the well-rounded retract contains interesting $\SL_n\BZ$-invariant subsets of dimension $\frac{n(n-1)}2$. For instance, recall that an element $A\in\SL_n\BR$ is well-rounded if the set of shortest non-zero vectors of the lattice $A\BZ^n$ generate $\BR^n$ as a vector space; equivalently, they generate, as a group, a finite index lattice of $A\BZ^n$. We will say that $A\in\SL_n\BR$ is {\em extremely well-rounded} if the shortest non-zero vectors of $A\BZ^n$ generate the whole lattice $A\BZ^n$. The subset $\CX'$ of $S_n$ consisting of extremely well-rounded elements is $\SL_n\BZ$-invariant and has dimension $\frac{n(n-1)}2$. While $\CX'=\CX$ for $n=2,3$ and $4$ the set $\CX'$ is a proper subset of the well-rounded retract for $n\ge 5$. In \cite{weird} we proved that $\CX'$ is not contractible for $n\ge 5$. This result follows now directly from Theorem \ref{main}:

\begin{kor} \cite{weird}\label{weird}
The subset $\CX'\subset S_n$ of extremely well-rounded elements is not contractible.\qed
\end{kor}

In order to prove Theorem \ref{main} it suffices to show that whenever $\CY$ is a closed proper $\SL_n\BZ$-invariant subset of $\CX$, there is a torsion-free, finite index subgroup $\Gamma\subset\SL_n\BZ$ such that the inclusion $\CY/\Gamma\hookrightarrow\CX/\Gamma$ is not a homotopy equivalence. We proceed as follows: First we show that there is $A\in\CX\setminus\CY$ with the property that there is a torsion-free, finite index subgroup $\Gamma$ of $\SL_n\BZ$ and a non-trivial homology class $[\alpha]\in H_{n-1}(\bar M_\Gamma,\D\bar M_\Gamma)$ represented by a cycle $\alpha$ which intersects the well-rounded retract exactly at $A$. Here $\bar M_\Gamma$ is the Borel-Serre compactification of the locally symmetric space $M_\Gamma=S_n/\Gamma$ and the homology is with coefficients in the ring $\BZ/2\BZ$. The class $[\alpha]$ is dual to some class $[\beta]\in H_{\frac{n(n-1)}2}(M_\Gamma)$. The fact that the cycle $\alpha$ does not intersect $\CY$ implies that $[\beta]$ is not in the image of $H_*(\CY/\Gamma)$ in $H_*(\CX/\Gamma)$. This shows that the inclusion $\CY/\Gamma$ in $\CX/\Gamma$ is not a homotopy equivalence.
\medskip

In \cite{weird}, we used this strategy to prove Corollary \ref{weird}. In that particular case we faced much simpler technical problems since it was possible to explicitly find a rational maximal flat intersecting $\CX$ exactly once, at a point outside of $\CX'$. Even in the case $n=2$, it is easy to see that for a generic point $A\in\CX$, every maximal flat through $A$ intersects $\CX$ many times. To bypass this problem we give an elementary, though somewhat involved, construction of the cycle $\alpha$.
\medskip

The paper is organized as follows: In Section \ref{sec:ss} we review some facts about the symmetric space $S_n=\SO_n\bs\SL_n\BR$ and its quotients. In Section \ref{sec:wr} we discuss some properties of the well-rounded retract, proving that a generic well-rounded element in $S_n$ has exactly $2n$ shortest vectors. In Section \ref{sec:homology} we show that certain homology classes are non-trivial; all the results in this section are surely well-known. In Section \ref{sec:disk} we derive Theorem \ref{main} from a result, Proposition \ref{meat}, proved in Section \ref{sec:systoles}. Proposition \ref{meat}, the key point of this paper, yields nontrivial cycles in $C_{n-1}(\bar M_\Gamma,\D\bar M_\Gamma)$ which intersect the well-rounded retract at a single point.
\medskip

We thank Mladen Bestvina for enduring us while we were working on this project and for suggesting the strategy behind the proof of Theorem \ref{main}. The first author is thankful to the University of Southampton and to Juan's mama for their hospitality during the realization of this paper. The second author thanks the Department of Mathematics of Stanford University. The results of this paper were obtained while the second author was a member of the Department of Mathematics of the University of Chicago. Finally, both authors thank the {\em scientific} organizers of the Fall 2007 programs {\em Geometric Group Theory} and {\em Teichm\"uller theory and Kleinian groups} at the MSRI.

%
%
%

\medskip

\noindent{\bf Notation.}
\medskip

We denote by $\{e_1,\dots,e_n\}$ and $\vert\cdot\vert$ the standard basis and euclidean norm of $\BR^n$. Sometimes we will write elements in $\BR^n$ as columns and sometimes as rows; we hope that this does not cause any confusion. If $U$ is a linear subspace of $\BR^n$, denote by $U^\bot$ its orthogonal complement with respect to the standard euclidean product. We will use the same symbol to denote both an equivalence class and a representative of the equivalence class. For example, we use the same notation for an element in $\SL_n\BR$ and for the corresponding element in the symmetric space $S_n=\SO_n\bs\SL_n\BR$, or in even smaller quotients such as  $S_n/\SL_n\BZ$. We will however consistently denote the homology class corresponding to a cycle $\alpha$ by $[\alpha]$. All the homology groups considered below have coefficients in the field $\BZ/2\BZ$ of two elements, although  everything remains true with respect to any other commutative ring with unit.

%
%
%

\section{The symmetric space $S_n=\SO_n\bs\SL_n\BR$}\label{sec:ss}
Up to scaling, the manifold $S_n=\SO_n\bs\SL_n\BR$ admits a unique symmetric metric invariant under the right action of $\SL_n\BR$; we shall always assume $S_n$ to be endowed with such a metric. 
The restriction of the right action of $\SL_n\BR$ on $S_n$ to $\SL_n\BZ$ is discrete. Moreover, any torsion-free subgroup $\Gamma$ of $\SL_n\BZ$ acts freely and hence the quotient $M_\Gamma=S_n/\Gamma$ is a smooth locally symmetric manifold. It is well-known that $\SL_n\BZ$ contains torsion-free finite index subgroups. If $\Gamma\subset\SL_n\BZ$ is any such subgroup, then the manifold $M_\Gamma$ is not compact, but is homeomorphic to  the interior of a compact manifold $\bar M_\Gamma$, the so-called Borel-Serre compactification of $M_\Gamma$ \cite{Borel-Serre}.

For every $v\in\BR^n$, the {\em length function}
\begin{equation}\label{eq:length-function}
l_v:S_n\to\BR,\ \ l_v(A)=\vert Av\vert
\end{equation}
is well-defined, analytic and convex. In particular we have
\begin{equation}\label{eq:convex}
l_v(A'')\le\max\{l_v(A),l_v(A')\}
\end{equation}
for all $A,A'\in S_n$ and every $A''$ in the unique geodesic segment $[A,A']$ joining $A$ and $A'$ in $S_n$. It should be observed that for every $B\in\SL_n\BR$ we have $l_v(AB)=l_{Bv}(A)$. Since $\SL_n\BZ$ acts on the set $\BZ^n\setminus\{0\}$, this implies that the function
\begin{equation}\label{eq:systole}
\syst_1:S_n\to(0,\infty),\ \ \syst_1(A)=\min_{v\in\BZ^n,v\neq 0}l_v(A)
\end{equation}
is $\SL_n\BZ$-invariant. The quantity $\syst_1(A)$ is said to be the {\em systole}, or {\em first minimum}, of $A\in S_n$. The elements of the set
\begin{equation}\label{eq:systoles}
\CS_1(A)=\{v\in\BZ^n \ \vert \ l_v(A)=\syst_1(A)\}
\end{equation}
are said to be the {\em systoles} or {\em shortest vectors} of $A$. 

Ash proved in \cite{Ash-Morse} that the systole function is a topological Morse function (see also Bavard \cite{Bavard} and Akrout \cite{Akrout}). Moreover, the induced function on $S_n/\SL_n\BZ$ is proper by Mahler's compactness theorem.

\begin{named}{Mahler's compactness theorem}
A closed subset $K\subset S_n/\SL_n\BZ$ is compact if and only if there is $\epsilon>0$ with
$\syst_1(A)\ge\epsilon$ for all $A\in K$.
\end{named}

We deduce from \eqref{eq:convex} and Mahler's compactness theorem the following important observation:

\begin{lem}\label{homotopy}
Let $\Gamma$ be a torsion-free subgroup of $\SL_n\BZ$, $N$ a manifold, and $f,g:N\to S_n$ two continuous maps such that for all $\epsilon>0$ there is a compact set $K_\epsilon\subset N$ with the following property:
\begin{itemize}
\item[(*)] For all $x\notin K_\epsilon$ there is $v\in\BZ^n\setminus\{0\}$ with $l_v(f(x)),l_v(g(x))<\epsilon$. 
\end{itemize}
Then the compositions of $f$ and $g$ with the projection $\pi:S_n\to M_\Gamma$ are properly homotopic.
\end{lem}
\begin{proof}
Let $H:N\times[0,1]\to S_n$ be the geodesic homotopy from $f$ to $g$, i.e. $t\to H_t(x)$ traverses with constant velocity the geodesic segment $[f(x),g(x)]$. We claim that $h=\pi\circ H$ is proper. Let $C$ be a compact subset of $M_\Gamma=S_n/\Gamma$. By Mahler's compactness theorem there is some $\epsilon$ positive with $\syst_1(A)\ge\epsilon$ for all $A\in C$. For such an $\epsilon$, let $K_\epsilon\subset N$ be the compact subset provided by (*). Then for $x\notin K_\epsilon$ there is some $v_x\in\BZ$, $v_x\neq 0$, with $l_{v_x}(f(x)),l_{v_x}(g(x))<\epsilon$. By \eqref{eq:convex} we have then $l_{v_x}(H_t(x))<\epsilon$ for all $t\in[0,1]$. This implies that $h^{-1}(C)\subset K_\epsilon\times[0,1]$, proving that it is proper.
\end{proof}

We will use Lemma \ref{homotopy} several times in the following situation. 

\begin{kor}\label{kor:homotopy}
Assume that $\Gamma$ is a finite index subgroup of $\SL_n\BZ$, and that $N\subset\SL_n\BR$ projects properly to $M_\Gamma=\SO_n\bs\SL_n\BR/\Gamma$. Then for every $B\in\SL_n\BR$ the projections of $N$ and of $BN=\{Bx, x\in N\}$ to $M_\Gamma$ are properly homotopic.\qed
\end{kor}

%
%
%

\section{The well-rounded retract}\label{sec:wr}
In this section we discuss briefly some of the properties of the well-rounded retract. Recall the definition of the systole \eqref{eq:systole} and of the set of systoles \eqref{eq:systoles} of a point  $A\in S_n$. Let also 
\begin{equation}\label{eq:span-syst}
\Lambda_1(A)=\Span_\BR(\CS_1(A))
\end{equation}
be the linear subspace of $\BR^n$ generated by the set of systoles of $A$.

\begin{defi*}
An element $A\in S_n$ is {\em well-rounded} if $\Lambda_1(A)=\BR^n$. The subset $\CX$ of $S_n$ consisting of all well-rounded elements is called the {\em well-rounded retract}.
\end{defi*}

As mentioned in the introduction, Soul\'e \cite{Soule} and Ash \cite{Ash} proved that $\CX$ is an $\SL_n\BZ$-invariant deformation retract. The idea behind this result is simple and beautiful, and so we explain it briefly here:

\begin{sat}[Soul\'e, Ash]\label{well-rounded}
The well-rounded retract $\CX$ is a deformation retract of $S_n$.
\end{sat}

For $k=1,\dots,n$ let $\CX_k$ be the set of those $A\in S_n$ for which we have $\dim\Lambda_1(A)\ge k$. We have the following chain of nested $\SL_n\BZ$-invariant subspaces:
$$\CX=\CX_n\subset\CX_{n-1}\subset\dots\subset\CX_1=S_n$$
In order to prove Theorem \ref{well-rounded} it suffices to show that for $k=1,\dots,n-1$ the space $\CX_{k+1}$ is an $\SL_n\BZ$-equivariant spine of $\CX_k$; we construct a retraction. Given $A\in\CX_k$ and $\lambda\in\BR$, consider the one-parameter family of linear maps
\begin{equation}\label{flow}
T_A^\lambda\in\SL_n\BR,\ \ \ T_A^\lambda(v)=\left\{
\begin{array}{ll}
e^{(n-k)\lambda}v  & \hbox{for}\ v\in A\Lambda_1(A)  \\
e^{-k\lambda}v  &  \hbox{for}\ v\in (A\Lambda_1(A))^\bot
\end{array}
 \right.
\end{equation}
In other words, for positive $\lambda$ the map $T_A^\lambda$ expands the subspace generated by the image of the shortest vectors of $A$, while contracting the orthogonal complement. Observe that for $U\in\SO_n$ we have $T_{UA}^\lambda UA=UT_A^\lambda A$; hence the point $T_A^\lambda A\in S_n$ depends only on $A$ and not on the choice of representative. 

Now $T_A^0 A=A$, and if $A\in\CX_k\setminus\CX_{k+1}$, there is some $\lambda$ positive with $T_A^\lambda A\in\CX_{k+1}$. For $A\in\CX_k$, let $\tau(A)\ge 0$ be maximal such that
$$T_A^\lambda A\in\CX_k\setminus\CX_{k+1}\ \ \hbox{for all}\ \lambda\in[0,\tau(A))$$
By definition $\tau(A)=0$ for $A\in\CX_{k+1}$. The function $A\mapsto\tau(A)$ is continuous on $\CX_k$, which implies that 
\begin{equation}\label{eq:homotopy}
[0,1]\times\CX_k\to\CX_k,\ \ (t,A)\mapsto T_A^{t\tau(A)}A
\end{equation}
is continuous as well. By definition, this homotopy is $\SL_n\BZ$-equivariant, starts with the identity, and ends with a projection of $\CX_k$ to $\CX_{k+1}$. This proves that $\CX_{k+1}$ is an $\SL_n\BZ$-equivariant spine of $\CX_k$ for $k=1,\dots,n-1$, concluding the sketch of the proof of Theorem \ref{well-rounded}.\qed
\medskip

It is not difficult to prove that $\CX_k$ is a co-dimension $k-1$ semi-algebraic set, i.e., that it is given by a locally finite collection of inequalities and (quadratic) algebraic equations. Hence $\CX$ is homeomorphic to a CW-complex of dimension $\dim(\CX)=\dim S_n-(n-1)=\frac{n(n-1)}2$. It is also easy to see that $\CX/\Gamma$ is compact. We prove now that a generic point in $\CX$ has exactly $2n$ shortest vectors:
 
\begin{prop}\label{prop:well-rounded}
The set of those $A\in\CX$ for which there are $v_1,\dots,v_n\in\BZ^n$ linearly independent with $\CS_1(A)=\{\pm v_1,\dots,\pm v_n\}$ is dense in $\CX$.
\end{prop}

In order to prove Proposition \ref{prop:well-rounded} we will use the following not very surprising but also not completely obvious geometric lemma.

\begin{lem}\label{fact}
Assume that $\CS$ is a finite subset of the sphere $\BS^{n-1}$ in $\BR^n$ with the property that $\BR^n=\Span_\BR\CS$ and assume that if $v\in\CS$ then $-v\in\CS$ as well. Then there is basis $\CB$ of $\BR^n$ contained in $\CS$ and a linear map $F:\BR^n\to\BR^n$ close to the identity such that for $v\in\CS$ we have $\vert Fv\vert=\vert v\vert$ if $\pm v\in\CB$ and $\vert Fv\vert>\vert v\vert$ otherwise.
\end{lem}

Assuming Lemma \ref{fact}, we prove Proposition \ref{prop:well-rounded}. Given $A\in\CX$ choose a representative in $\SL_n\BR$, again denoted by $A$. By definition, the image $A\CS_1(A)$ of the set of systoles of $A$ generates $\BR^n$ and is contained in the round sphere $\BS^{n-1}_{\syst_1(A)}$ of radius $\syst_1(A)$. Let $\CB\subset A\CS_1(A)$ and $F:\BR^n\to\BR^n$ be the basis and the linear map provided by Lemma \ref{fact}. We set $A^{-1}\CB=\{v_1,\dots,v_n\}$ and $A'=\frac 1{\sqrt[n]{\det(F)}}FA$. Since we may assume that $F$ is very close to the identity, we have that $A'$ is very close to $A$, and hence $\CS_1(A')\subset\CS_1(A)$. It follows now from Lemma \ref{fact} that $\CS_1(A')=\{\pm v_1,\dots,\pm v_n\}$. This concludes the proof of Proposition \ref{prop:well-rounded}. \qed 
\medskip

We prove now Lemma \ref{fact}:

\begin{proof}[Proof of Lemma \ref{fact}]
We will prove the lemma by induction on the number of elements in $\CS$. There is nothing to show if $\CS$ has $2n$ elements. Assume that we have proved it for all sets with at most $2k\ge 2n$ elements and that $\CS$ has $2(k+1)$ elements. We begin observing that there is a co-dimension one linear subspace $U\subset\BR^n$ generated by $U\cap\CS$ and such that there are at least four elements in $\CS$ which don't belong to $U$; recall that whenever $v\in\CS$ then $-v\in\CS$ as well. 

We choose now $v\in\CS$, $v\notin U$ with minimal angle $\angle(U,v)=\theta\in(0,\frac\pi 2)$. Let $V$ be the codimension one linear subspace containing $v$ and the intersection $(\BR v)^\perp\cap U$ of the orthogonal complement of $\BR v$ and $U$. The choice of $v$ and the construction of $V$ imply that $\pm v$ are the only elements in $\CS\cap V$ which don't belong to $U$. The planes $U$ and $V$ have angle $\theta$ and divide $\BR^n$ into two open sectors, $C_1$ and $C_2$,  with angle $\theta$ and two also open sectors, $C_3$ and $C_4$, with angle $\pi - \theta$. Moreover we have $\CS\cap(C_1\cup C_2)=\emptyset$ but $\CS\cap(C_3\cup C_4)\neq\emptyset$, by the minimality of $\theta$.

For $\eta>\theta$ with $\eta-\theta$ small we can consider the linear map $F:\BR^n\to\BR^n$ which is the identity on $U$, an isometry when restricted to $V$, and which opens $C_1$ and $C_2$ to angle $\eta$. The map $F$ preserves the length of points in $U\cup V$, reduces the length of vectors in $C_1\cup C_2$ and increases the length of vectors in $C_3\cup C_4$. In particular, $F$ maps $\CS\cap(U\cup V)$ to a subset of $\BS^n$ which still generates $\BR^n$ and increases the length of the (at least two) remaining points in $\CS$. This concludes the induction step and the proof of Lemma \ref{fact}.\end{proof}

%
%
%

\section{A bit of homology}\label{sec:homology}
In this section we give elementary proofs of some homological results which are probably well-known to experts and non-experts alike.
\medskip

As mentioned above, $\SL_n\BZ$ contains torsion-free subgroups of finite index, and any such subgroup acts freely and discretely on $S_n$; as always, we denote the quotient manifold by $M_\Gamma=S_n/\Gamma$ and its Borel-Serre compactification by $\bar M_\Gamma$. If $U\subset\bar M_\Gamma$ is a regular neighborhood of $\D\bar M_\Gamma$, we have $H_*(\bar M_\Gamma,U)\simeq H_*(\bar M_\Gamma,\D\bar M_\Gamma)$. In particular, we can consider every properly immersed submanifold of $M_\Gamma$ as a cycle in $C_*(\bar M_\Gamma,\D\bar M_\Gamma)$. Recall that we always consider homology with coefficients in $\BZ/2\BZ$. 

Before stating the main result of this section, we recall that by Lefschetz duality there is a non-degenerate pairing
$$\iota:H_{n-1}(\bar M_\Gamma,\D\bar M_\Gamma)\times H_{\frac{n(n-1)}2}(M_\Gamma)\to\BZ/2\BZ$$
which can be computed as follows. Given homology classes $[\alpha]\in H_{n-1}(\bar M_\Gamma,\D\bar M_\Gamma)$ and $[\beta]\in H_{\frac{n(n-1)}2}(M_\Gamma)$, represent them by cycles $\alpha$ and $\beta$ in general position. Then $\iota([\alpha],[\beta])$ is just the parity of the cardinality of the set $\alpha\cap\beta$. Observe that in order to prove that a cycle $\beta\in C_{\frac{n(n-1)}2}(M_\Gamma)$ represents a non-trivial homology class, it suffices to find a cycle $\alpha\in C_{n-1}(\bar M_\Gamma,\D\bar M_\Gamma)$ which intersects $\beta$ transversally at a single point; if this is the case we will say that the two classes $[\alpha]$ and $[\beta]$ are dual to each other. This is the argument used in \cite{weird} to prove:

\begin{prop}\label{homology-weird}
Let $\Gamma$ be a finite index torsion-free subgroup of $\SL_n\BZ$, $\Delta$ the connected component of the identity in the diagonal subgroup of $\SL_n\BR$ and $\Nil$ the subgroup of $\SL_n\BR$ consisting of upper triangular matrices with units in the diagonal. Then the projection of $\Delta$ and $\Nil$ to $M_\Gamma$ represent dual, and hence nontrivial, homology classes in $H_{n-1}(\bar M_\Gamma,\D\bar M_\Gamma)$ and $H_{\frac{n(n-1)}2}(M_\Gamma)$, respectively.\qed
\end{prop}

Proposition \ref{homology-weird} is surely well-known, as is the following slightly more general version.

\begin{kor}\label{homology}
Given $B\in\GL_n\BQ$ assume that $\Gamma\subset\SL_n\BZ$ is a finite index torsion-free subgroup with $B^{-1}\Gamma B\subset\SL_n\BZ$, and that $\Delta$ and $\Nil$ are as in Proposition \ref{homology-weird}. Then the projections of $B\Delta B^{-1}$ and $B\Nil B^{-1}$ to $M_\Gamma$ represent dual, and hence nontrivial, homology classes in $H_{n-1}(\bar M_\Gamma,\D\bar M_\Gamma)$ and $H_{\frac{n(n-1)}2}(M_\Gamma)$, respectively.
\end{kor}
\begin{proof}
The map $\phi:S_n\to S_n$, $\phi(X)=XB^{-1}$ induces a diffeomorphism $\Phi:M_{B^{-1}\Gamma B}\to M_\Gamma$. By Proposition \ref{homology-weird} the projections of $\Delta$ and $\Nil$ represent dual homology classes in $M_{B^{-1}\Gamma B}$. Pushing forward with $\Phi$, we obtain dual cycles $\Delta B^{-1}$ and $\Nil B^{-1}$. By Corollary \ref{kor:homotopy}, these cycles are properly homotopic, and hence homologous, to the cycles $B\Delta B^{-1}$ and $B\Nil B^{-1}$. The claim follows.
\end{proof}

%
%
%

\section{Proof of Theorem \ref{main}}\label{sec:disk}
In the next section we will show:

\begin{prop}\label{meat}
Assume that $A\in\CX$ is such that there are $v_1,\dots,v_n\in\BZ^n$ linearly independent with $\CS_1(A)=\{\pm v_1,\dots,\pm v_n\}$. Let $B\in\GL_n\BQ$ be the matrix with columns $v_1,\dots,v_n$, and let $\Gamma$ be a finite index torsion-free subgroup of $\SL_n\BZ\cap B\SL_n\BZ B^{-1}$. Then the non-trivial homology class $[B\Delta B^{-1}]$ is represented by a cycle $\alpha\in C_{n-1}(\bar M_\Gamma,\D\bar M_\Gamma)$ whose support intersects the well-rounded retract $\CX$ only in $A$.
\end{prop}

Assuming Proposition \ref{meat}, we prove Theorem \ref{main}:

\begin{named}{Theorem \ref{main}}
The well-rounded retract $\CX$ is a minimal $\SL_n\BZ$-invariant spine of the symmetric space $S_n=\SO_n\bs\SL_n\BR$.
\end{named}

\begin{proof}
Assume that $\CY\subset\CX$ is a proper, closed, $\SL_n\BZ$-invariant subset of $\CX$. As mentioned in the introduction, in order to show that $\CY$ is not contractible, it suffices to prove that for some $\Gamma\subset\SL_n\BZ$ the induced map $\CY/\Gamma\to\CX/\Gamma$ is not a homotopy equivalence.

By Proposition \ref{prop:well-rounded} there is $A\in\CX\setminus\CY$ and a linearly independent subset $\{v_1,\dots,v_n\}\subset\BZ$ with $\CS_1(A)=\{\pm v_1,\dots,\pm v_n\}$. Let $B\in\GL_n\BQ$ be the matrix with columns $v_1,\dots,v_n$. The subgroups $\SL_n\BZ$ and $B\SL_n\BZ B^{-1}$ are commensurable and hence there is a torsion-free finite index subgroup $\Gamma\subset\SL_n\BZ\cap B\SL_n\BZ B^{-1}$. By Proposition \ref{meat}, the homology class $[B\Delta B^{-1}]\in H_{n-1}(\bar M_\Gamma,\D\bar M_\Gamma)$ is represented by a cycle $\alpha$ with $\alpha\cap\CX=\{A\}$. On the other hand, the class $[B\Delta B^{-1}]$ is dual to some class $[\beta]\in H_{\frac{n(n-1)}2}(M_\Gamma)$ by Corollary \ref{homology}. Since $\alpha$ represents $[B\Delta B^{-1}]$ and intersects $\CX$ only at $A$, we deduce that every cycle contained in $\CX/\Gamma$ and representing $[\beta]$ has to contain $A$ in its support. In particular, the map 
$$H_{\frac{n(n-1)}2}(\CY/\Gamma)\to H_{\frac{n(n-1)}2}(\CX/\Gamma)$$
is not surjective. This implies that the map $\CY/\Gamma\to\CX/\Gamma$ is not a homotopy equivalence.
\end{proof}

%
%
%

\section{Flags of systoles}\label{sec:systoles}
In this section we prove Proposition \ref{meat}. The first step is to construct a certain continuous map 
\begin{equation}\label{eq:map0}
\Phi:S_n\times[0,\infty)\to S_n
\end{equation}
which essentially pushes points in $S_n\setminus\CX$ away from $\CX$. 

To begin with, recall the definition of the systole $\syst_1(A)$ of $A\in S_n$. We can extend this definition as follows: for $i=1,\dots,n$, the {\em i-th systole} of $A$ is given by
\begin{equation}\label{eq:i-systole}
\syst_i(A)=\inf\{r\vert \dim_\BR(\Span_\BR\{v\in\BZ\ \hbox{with}\ \vert Av\vert< r\})\ge i\}
\end{equation}
In other words, $\syst_i(A)$ is the infimum of those $r$ for which the set of vectors $v$ in $\BZ^n$ whose image $Av$ has length less than $r$ generates an $i$-dimensional subspace of $\BR^n$. Equivalently,
\begin{equation}\label{eq:i-systole2}
\syst_i(A)=\sup\{r\vert \dim_\BR(\Span_\BR\{v\in\BZ\ \hbox{with}\ \vert Av\vert< r\})< i\}
\end{equation}
The $i$-th systole coincides with Minkowski's $i$-th successive minimum of the lattice $A\BZ^n$ with respect to the ball $B_1$ of radius $1$ in $\BR^n$. See \cite{Martinet} for more about successive minima.

For $i=1,\dots,n$, the $i$-th systole function
$$\syst_i:S_n\to(0,\infty)$$
is well-defined and $\SL_n\BZ$-equivariant. We claim that it is continuous. In fact, if $(A_k)$ is a sequence in $S_n$ converging to some $A\in S_n$ then for all $r$ the finite sets $\{v\in\BZ^n,\vert A_kv\vert<r\}$ converge in the Gromov-Hausdorff topology to the (again finite) set $\{v\in\BZ^n,\vert Av\vert<r\}$. Since $\BZ^n$ is discrete, we have that for all sufficiently large $k$
$$\{v\in\BZ^n,\vert A_kv\vert<r\}=\{v\in\BZ^n,\vert Av\vert<r\}$$
Together with (\ref{eq:i-systole}), this implies that $\syst_i$ is lower semi-continuous. Likewise \eqref{eq:i-systole2} and the same argument yield upper semi-continuity.

\begin{lem}\label{i-syst-cont}
The function $\syst_i:S_n\to(0,\infty)$ is continuous and $\SL_n\BZ$-equivariant for $i=1,\dots,n$.\qed 
\end{lem}

Recall now the definition of $\Lambda_1(A)$ given in \eqref{eq:span-syst}. We extend this definition,  setting for $i=1,\dots,n$
$$\Lambda_i(A)=\Span_\BR(\{v\in\BZ^n,\vert Av\vert\le\syst_i(A)\})$$
In order to avoid treating special cases we set $\Lambda_0(A)=0$ for all $A\in S_n$. By definition 
\begin{equation}\label{flag}
0\subsetneq\Lambda_1(A)\subset\dots\subset\Lambda_n=\BR^n
\end{equation}
and $\dim_\BR(\Lambda_i(A))\ge i$. Observe that for $i<n$ this last inequality is strict if $A$ is well-rounded. In particular, we cannot expect that the subspaces $\Lambda_i(A)$ depend continuously of $A$. However we have the following weak continuity, which can be proved with essentially the same argument as Lemma \ref{i-syst-cont}:

\begin{lem}\label{flag-cont}
Assume that $(A_k)$ is a sequence in $S_n$ converging to some $A\in S_n$. Then there is $k_0$ such that for all $k\ge k_0$ and $i\in\{1,\dots,n\}$ there is a unique $\kappa(k,i)\in\{1,\dots,n\}$ with 
\begin{itemize}
\item $\Lambda_{\kappa(k,i)}(A_k)=\Lambda_i(A)$, and 
\item if $\kappa(k,i)\neq n$ then $\Lambda_{\kappa(k,i)+1}(A_k)\neq\Lambda_i(A)$. 
\end{itemize}
If moreover $i'$ is minimal with $\syst_{i'}(A)=\syst_i(A)$ then 
$$\lim_{k\to\infty}\syst_{j_k}(A_k)=\syst_i(A)$$
for all choices of $j_k$ with $\kappa(k,i'-1)<j_k\le\kappa(k,i)$.\qed
\end{lem}

We use the flag \eqref{flag} to construct the continuous map \eqref{eq:map0}. To begin with we consider for $i=1,\dots,n$ the subspace
$$\Theta_i(A)=(A\Lambda_{i-1}(A))^\perp\cap(A\Lambda_i(A))$$
In more plain language, $\Theta_i(A)$ is the orthogonal complement of the image of $\Lambda_{i-1}(A)$ under $A$ within the image of $\Lambda_i(A)$. We have thus the orthogonal decomposition 
\begin{equation}\label{orthogonal}
\BR^n=\Theta_1(A)\oplus\dots\oplus\Theta_n(A)
\end{equation}
together with the associated orthogonal projections
\begin{equation}\label{eq:projection}
\pi_{\Theta_i(A)}:\BR^n\to\Theta_i(A)
\end{equation}
We define now for $x\in\BR^n$
\begin{equation}\label{eq:map}
\Phi_t(A)x=\frac 1{\sqrt[n]{\prod_{i=1}^n\syst_i(A)^{t\dim_\BR\Theta_i(A)}}}\sum_{i=1}^n\syst_i(A)^t\pi_{\Theta_i(A)}(Ax)
\end{equation}
The multiplicative factor in \eqref{eq:map} ensures that $\Phi_t(A)\in\SL_n\BR$ for all $A\in\SL_n\BR$. Moreover, for all $U\in\SO_n$ we have $\Phi_t(UA)=U\Phi_t(A)$. In particular, we have a well-defined map
\begin{equation}\label{eq:map2}
\Phi_t:S_n\times[1,\infty)\to S_n
\end{equation}
It is easy to check that the map \eqref{eq:map2} is $\SL_n\BZ$-equivariant, and its continuity follows from Lemma \ref{flag-cont}. Moreover, since $\syst_1(A)\le\syst_i(A)$ for all $i$, we have for all $x\in\BR^n$
\begin{equation}\label{eq:bound}
\vert\Phi_t(A)x\vert\ge\left(\frac {\syst_1(A)}{\sqrt[n]{\prod_{i=1}^n\syst_i(A)^{\dim_\BR\Theta_i(A)}}}\right)^t\vert Ax\vert
\end{equation}
with equality if and only if $x\in\Lambda_1(A)$. In particular we see that $\Lambda_1(\Phi_t(A))=\Lambda_1(A)$ for all $t\ge 0$. Moreover, if $\Lambda_1(A)\neq\BR^n$ then the exponentiated quantity in \eqref{eq:bound} is less than $1$ and hence
$$\lim_{t\to\infty}\syst_1(\Phi_t(A))=0$$
On the other hand, if $\Lambda_1(A)=\BR^n$ then $\Phi_t(A)=A$ for all $t$.

Summing up, we have:

\begin{prop}\label{map}
There is a continuous map $\Phi:S_n\times[0,\infty)\to S_n$, $\Phi(A,t)=\Phi_t(A)$, with the following properties:
\begin{itemize}
\item $\Phi_0(\cdot)=\Id$,
\item $\Phi_t(A)\in\CX$ if and only if $A\in\CX$, and
\item if $A\notin\CX$ then $\lim_{t\to\infty}\vert\Phi_t(A)v\vert=0$ for all $v\in\Lambda_1(A)$.\qed
\end{itemize}
\end{prop}

We are now ready to prove Proposition \ref{meat}:

\begin{named}{Proposition \ref{meat}}
Assume that $A\in\CX$ is such that there are $v_1,\dots,v_n\in\BZ^n$ linearly independent with $\CS_1(A)=\{\pm v_1,\dots,\pm v_n\}$, let $B\in\GL_n\BQ$ be the matrix with columns $v_1,\dots,v_n$ and $\Gamma$ a finite index torsion-free subgroup in $\SL_n\BZ\cap B\SL_n\BZ B^{-1}$. Then the non-trivial homology class $[B\Delta B^{-1}]$ is represented by a cycle $\alpha\in C_{n-1}(\bar M_\Gamma,\D\bar M_\Gamma)$ whose support intersects the well-rounded retract $\CX$ only at $A$.
\end{named}

Recall that $\Delta$ is the connected component of the identity in the diagonal subgroup of $\SL_n\BR$.

\begin{proof}
In order to construct the cycle $\alpha$ we start with the map
$$g_1:\Delta\to M_\Gamma,\ \ g_1(X)=BXB^{-1}$$
By Proposition \ref{homology}, the cycle $g_1(\Delta)$ represents a non-trivial homology class in $H_{n-1}(\bar M_\Gamma,\D\bar M_\Gamma)$. The point $A$ may not belong to the image of $g_1(\Delta)$, but this can be easily corrected by considering the map
$$g_2:\Delta\to M_\Gamma,\ \ g_2(X)=ABXB^{-1}$$
Corollary \ref{kor:homotopy} implies that $g_1(\Delta)$ and $g_2(\Delta)$ are properly homotopic and hence homologous.

Now we have $g_2(\Id)=A$, but it is not clear at all how many other times  $g_2(\Delta)$ may intersect $\CX$. We correct this problem by constructing a third map $g_3$ properly homotopic to $g_2$. Before going further we identify $\Delta$ with $\BR^{n-1}$ via the following map
\begin{equation}\label{eq:diagonal}
(a_1,\dots,a_{n-1})\mapsto
\left(
\begin{array}{ccccc}
e^{a_1}  & 0  & \dots & 0 & 0 \\
0  & e^{a_2}  & \dots & 0 & 0 \\
\vdots  & \vdots  & \ddots & \vdots & \vdots \\
0 & 0 & \dots & e^{a_{n-1}} & 0 \\
0 & 0 & \dots & 0 & e^{-a_1-a_2-\dots-a_{n-1}} 
\end{array}
\right)
\end{equation}
A simple computation shows: 
\begin{lem}\label{lem:little-disk}
There is some $\epsilon>0$ such that for all $x\in B_\epsilon\subset\BR^{n-1}=\Delta$, \item $g_2(x)\in\CX$ if and only if $x=0$. If moreover $x\in B_\epsilon$, $x\neq 0$ and $v\in\CS_1(x)$ then we have
\begin{eqnarray}\label{limitg2}
\lim_{t\to\infty}l_v(g_2(tx))=0
\end{eqnarray}
Here $B_\epsilon$ is the ball of radius $\epsilon$ centered at $0$ in $\BR^{n-1}\simeq\Delta$.\qed
\end{lem}

We can now define the map $g_3:\BR^{n-1}\to M_\Gamma$. With $\epsilon$ as in Lemma \ref{lem:little-disk} and $\Phi$ the map provided by Proposition \ref{map}, we set
$$g_3(x)=
\left\{
\begin{array}{cc}
g_2(x)  & \vert x\vert\le\epsilon \\
\Phi_{\vert x\vert-\epsilon}(g_2(\frac x{\vert x\vert})) & \vert x\vert\ge\epsilon   
\end{array}
\right\}$$
In other words we extend radially, using the map $\Phi$ and the restriction of $g_2$ to $B_\epsilon$. Since $g_2(x)\notin\CX$ for $x$ with $\vert x\vert=\epsilon$, we deduce from Proposition \ref{map} that $g_3(x)\notin\CX$ for all $x$ with $\vert x\vert\ge\epsilon$. On the other hand, for $\vert x\vert\le\epsilon$ we have $g_3(x)=g_2(x)$. Hence
$$g_3(\BR^{n-1})\cap\CX=\{A\}$$
If $v\in\BZ^n$ is a systole for $g_2(x)$ with $\vert x\vert=\epsilon$, then we have by  
(\ref{limitg2})
$$\lim_{t\to\infty}l_v(g_2(tx))=0$$
and by Proposition \ref{map}
$$\lim_{t\to\infty}l_v(g_3(tx))=\lim_{t\to\infty}l_v(\Phi_{t-1}(g_2(x))=0$$
Lemma \ref{homotopy} implies now that the maps $g_2$ and $g_3$ are properly homotopic to each other. Hence the cycle $\alpha=g_3(\Delta)$ represents the non-trivial homology class $[B\Delta B^{-1}]\in H_{n-1}(\bar M_\Gamma,\D\bar M_\Gamma)$ and $\alpha\cap\CX=\{A\}$.
\end{proof}

\vspace{.2cm}

{\sc \tiny \noindent
Alexandra Pettet, Department of Mathematics, Stanford University}

{\sc \tiny \noindent
Juan Souto, Department of Mathematics, University of Michigan, Ann Arbor}

\end{document}